\documentclass[a4paper,12pt]{article}
\usepackage{graphicx}
\usepackage{amssymb}
\usepackage{amsmath}
\usepackage{mathtools}
\usepackage{amsfonts}
\usepackage{appendix}
\usepackage{color}
\usepackage{multirow}
\usepackage{pdflscape}
\usepackage{xcolor}
\usepackage{setspace}
\newcommand{\1}{\boldsymbol{1}}

\newtheorem{theorem}{Theorem}

\newtheorem{lemma}[theorem]{Lemma}

\newenvironment{proof}[1][Proof]{\textbf{#1.} }{\ \rule{0.5em}{0.5em}}

\begin{document}

\title{\bf Developing multivariate distributions using Dirichlet generator}

\bigskip

\author{M. Arashi , A. Bekker, D. de Waal and S. Makgai\\
	Shahrood University of Technology and University of Pretoria
 }

\date{}
\maketitle

\begin{quotation}
\noindent {\it Abstract:}
There exist several endeavours proposing a new family of extended distributions using the beta-generating technique. This is a well-known	mechanism in developing flexible distributions, by embedding the cumulative distribution function (cdf) of a baseline distribution within the beta distribution that acts as a generator. Univariate beta-generated distributions	offer many fruitful and tractable properties, and have applications in hydrology, biology and environmental sciences amongst other fields. In the univariate cases, this extension works well, however, for multivariate cases	the beta distribution generator delivers complex expressions.
In this chapter the proposed extension from the univariate to the multivariate domain addresses the need of flexible multivariate distributions that can	model a wide range of multivariate data. This new family of multivariate	distributions, whose marginals are beta-generated distributed, is constructed with the function $H(x_{1},...,x_{p})=F\left(  G_{1}(x_{1}),G_{2}%
(x_{2}),...,G_{p}(x_{p})\right)  $, where $G_{i}(x_{i})$ are the cdfs of the gamma (baseline) distribution and $F(\cdot)$ as the cdf of the Dirichlet distribution. Hence as a main example, a general model having the support 	$[0,1]^{p}$ (for $p$ variates), using the Dirichlet as the generator, is 	developed together with some distributional properties, such as the moment generating function. The proposed Dirichlet-generated distributions can be	applied to compositional data. The parameters of the model are estimated by using the maximum likelihood method. The effectiveness and prominence of the	proposed family is illustrated through analyzing simulated as well as two real datasets. A new model testing technique is introduced to evaluate the performance of the multivariate models.
\par

\end{quotation}\par
\section{Introduction}
\label{sec:1}
In many of the problems of interest to scientists, data consists of proportions and thus are subject to non-negativity and unit-sum constraints. Examples of such data can be found when analyzing rock compositions, household budgets, pollution components to name a few. Datasets such as these are known as compositional datasets and arise naturally in a great variety of disciplines such as biology, medicine, chemistry, economics, psychology, environmetrics, psychology and many others. The most widely studied distribution on the simplex is the Dirichlet distribution \cite{Balakrishnan2003}. Various generalizations of the Dirichlet distribution are proposed in literature, for
example see \cite{Connor and Mosimaan}, \cite{Barndorff-Nielsen},
\cite{Ehlers}, \cite{Thomas}, \cite{Epaillard} and \cite{Favaro}. For an
extensive review see \cite{Ng} and \cite{Kotz}. In particular, the Liouville distribution has been widely studied (see \cite{Gupta and Richards}). Specifically, a flexible Dirichlet was proposed by \cite{Ongaro}, by extending the basis of gamma independent random variables which generates the Dirichlet distribution. The Dirichlet prior is widely used in estimating discrete distributions and functionals of discrete distributions, and in fact the Dirichlet distribution is the conjugate prior of the categorical distribution and multinomial distribution.

In this chapter we propose a general multivariate construction methodology
using the Dirichlet probability density function (pdf) as the generator. This Dirichlet-generated class serves as good alternatives to the Dirichlet and generalized Dirichlet distributions for the statistical representation of
specific proportional data. This class is an evolution from the univariate framework describes below into a multivariate setting:%

\begin{equation}
H(x)=\int_{0}^{G(x)}f(y)dy,\label{Mother}%
\end{equation}

\noindent with pdf%

\begin{equation}
h(x)=f\left(  G(x)\right)  g(x),\label{Mother pdf}%
\end{equation}
\noindent where $G(\cdot)$ is a continuous cumulative distribution function
(cdf) and $f(\cdot)$ is the pdf of a random variable with support $[0,1].$ By introducing extra parameters in $f(\cdot)$ and $G(\cdot)$ the resulting distribution provides greater flexibility in adapting modality and skewness. \cite{Eugene} was the first to introduce the family of beta-generated normal distribution with $f(y)=y^{\alpha-1}(1-y)^{\beta-1}/B\left(
\alpha,\beta\right)  $ as the pdf of the well-known beta distribution, where
$B\left(  \alpha,\beta\right)  =\Gamma\left(  \alpha\right)  \Gamma\left(
\beta\right)  /\Gamma\left(  \alpha+\beta\right)  $ denotes the classical beta
function and $\Gamma\left(  \alpha\right)  =%
{\displaystyle\int\limits_{0}^{\infty}}
v^{\alpha-1}e^{-v}dv$ is the gamma function defined for all $\alpha>0.$ The resulting cdf and pdf are respectively%

\begin{equation}
H(x)=\frac{1}{B(\alpha,\beta)}\int_{0}^{G(x)}y^{\alpha-1}(1-y)^{\beta
	-1}dy\label{beta normal cdf}%
\end{equation}

\noindent and%

\begin{equation}
h(x)=\frac{1}{B(\alpha,\beta)}g(x)G^{\alpha-1}(x)[1-G(x)]^{\beta
	-1},\label{beta normal pdf}%
\end{equation}

\noindent where $\alpha>0,\beta>0,$ and $g(\cdot)$ and $G(\cdot)$ are the pdf and cdf respectively. The beta distribution $f(\cdot)$ is referred to as the generator and $G(\cdot)$ as the baseline distribution. Another development of $\left(  \ref{beta normal pdf}\right)  $ is based on the $i$th order statistic in a random sample of $n$ from a distribution $G(\cdot)$ with pdf $\left\{  n!/[(i-1)!(n-i)!]\right\}  g(x)G^{i-1}(x)[1-G(x)]^{n-1}$ where \cite{Jones 2004} extended the pdf of the $i$th order statistic by allowing $a=i$ and $b=n+1-i$ which is the pdf in (\ref{beta normal pdf}). Note that the relation $X=G^{-1}(F(\cdot))$ with $F(\cdot)$ being a beta-distributed random variable, can be used to simulate $X$ values. It is clear that special choices of the baseline model $G(\cdot)$ yield specific models generated by the classic beta distribution. In recent years, several scholars have shown great interest in defining new generalized classes of univariate continuous distributions by using this \textquotedblleft mother technique\textquotedblright\ (see $\left(
\ref{Mother}\right)  $) to generate new models. The interested reader is
referred to \cite{Elgarhy 2016} (and the references therein), \cite{Makgai
	2017}, \cite{Alexander},\ \cite{Barreto-Souza}, \cite{Nadarajah2006},
\cite{Zografos}, \cite{Mameli} \ and \cite{Nassar} for related studies, amongst others.

Mimicking the same construction methodology (\ref{Mother}), three classes of extended bivariate distributions with the beta as generator, can be obtained as follows:

\begin{itemize}
	\item Builder 1:
\end{itemize}%

\begin{equation}
H(x_{1},x_{2})=\frac{1}{B(\alpha,\beta)}\int_{0}^{G(x_{1})G(x_{2})}%
y^{\alpha-1}(1-y)^{\beta-1}dy\label{Builder 1}%
\end{equation}

\begin{itemize}
	\item Builder 2:
\end{itemize}%

\begin{equation}
H(x_{1},x_{2})=\frac{1}{B(\alpha,\beta)}\int_{0}^{G_{1}(x_{1})G_{2}(x_{2}%
	)}y^{\alpha-1}(1-y)^{\beta-1}dy\label{Builder 2}%
\end{equation}

\begin{itemize}
	\item Builder 3:
\end{itemize}%

\begin{equation}
H(x_{1},x_{2})=\frac{1}{B(\alpha,\beta)}\int_{0}^{G^{\ast}(x_{1},x_{2}%
	)}y^{\alpha-1}(1-y)^{\beta-1}dy\label{Builder 3}%
\end{equation}

\noindent where $G_{i}(\cdot)$, $i=1,2$, can be any cdf of a baseline univariate distrbution and $G^{\ast}%
(\cdot,\cdot)$ is the cdf of the baseline bivariate distribution, $\alpha>0,\beta>0$.
\pagebreak

From Builder 1, the pdf has the form%
\begin{eqnarray}
h(x_{1},x_{2})&=&\frac{1}{B(\alpha,\beta)}G^{\alpha-1}(x_{1})G^{\alpha-1}
(x_{2})[1-G(x_{1})]^{\beta-1}[1-G(x_{2})]^{\beta-1}\\ \nonumber
&\times& [g(x_{1})G(x_{2}%
)+G(x_{1})g(x_{2})],\label{Builder 1 pdf}%
\end{eqnarray}

\noindent where $g(\cdot)$ is the pdf relative to the cdf $G(\cdot).$ In this case only one cdf contributes as baseline to develop the bivariate
distribution and is a special case of Builders 2 and 3. The advantage of
Builder 1 compared to Builder 2, is that it has fewer number of parameters. Makgai et al (2019) proposed Builder 3 and studied the properties and dependence structure of the class formed along with multivariate beta-generated distribution. Samanthi and Sepanski (2017) employed copulas to construct a bivariate extension of beta-generated distributions.

From completely a different viewpoint, \cite{Sarabia 2014} formed a bivariate distribution (see also \cite{Ristic}), using the \cite{Olkin}  beta pdf as generator:%

\begin{eqnarray}
h(x_{1},x_{2})&=&\frac{1}{B(\alpha,\beta,\gamma)}g_{1}(x_{1})g_{2}(x_{2}%
) \\ \nonumber
&\times& \frac{G_{1}^{\alpha-1}(x_{1})G_{2}^{\beta-1}(x_{2})[1-G_{1}(x_{1}%
	)]^{\beta+\gamma-1}[1-G_{2}(x_{2})]^{\alpha+\gamma-1}}{[1-G_{1}(x_{1}%
	)G_{2}(x_{2})]^{\alpha+\beta+\gamma}}.\label{Sarabia}%
\end{eqnarray}

However, the purpose of this study is not to study Builders 1-3, but to propose a general multivariate construction methodology using the Dirichlet pdf as the generator, with the baseline as the product of independent cdfs. This range of baseline distributions can be the exponential, Weibull, gamma, Fr\'{e}chet, etc. Suppose that $G(\cdot)$ belongs to the Pareto class, then $H(\cdot)$ is referred to as the Dirichlet-Pareto distribution function. The introduction of the Dirichlet distribution as the generating distribution $F(\cdot)$, creates the opportunity to apply a wide range of multivariate distributions. In this context, Section 2 provides the basic elements of the construction, that will be described in Section 3, with specific emphasis on the Dirichlet-Gamma distribution. In Section 4  some properties of the newly proposed multivariate distribution are  discussed. To illustrate the effectiveness of the latter model, the well-known Dirichlet distribution is compared to the Dirichlet-Gamma distribution via a simulation studies and an analysis of real datasets using different measures. Finally, some conclusions are given in Section 5.

\section{Ingredients}
\label{sec2}
In this section, the basic notation and definitions (ingredients) underlying the construction that will described in Section 3, are recalled.  A random vector $\boldsymbol{Y}=(Y_{1},\ldots,Y_{p})\in\mathcal{R}^{p}$ is said
to have Dirichlet distribution (or standard Dirichlet) with parameters
$\boldsymbol{\alpha}=$ $(\alpha_{1},\cdots,\alpha_{p};\alpha_{p+1})$ for
\ $\alpha_{i}>0,i=1,...,p+1,$ $p\geq2.$, if the pdf is given by%

\[
f(\boldsymbol{y})=\frac{\Gamma\left(  \alpha_{+}\right)  }{\Gamma(\alpha
	_{1})\cdots\Gamma(\alpha_{p+1})}y_{1}^{\alpha_{1}-1}\cdots y_{p}^{\alpha
	_{p}-1}\left(  1-\sum_{i=1}^{p}y_{i}\right)  ^{\alpha_{p+1}-1},%
\]

\noindent where $y_{i}>0,i=1,...,p,\sum_{i=1}^{p}y_{i}<1,$ use $\alpha
_{+}=\sum_{i=1}^{p+1}\alpha_{i}.$

For convenience, denote $Y_{p+1}=1-\sum_{i=1}^{p}Y_{i}, \boldsymbol{Y}%
^{^{\prime}}=(Y_{1},\ldots,Y_{p};Y_{p+1})=(\boldsymbol{Y};Y_{p+1})$ and write the above Dirichlet distribution as $\boldsymbol{Y\sim}Dir(\boldsymbol{\alpha
})\boldsymbol{,}$ or simply $\boldsymbol{Y}^{^{\prime}}\boldsymbol{\sim
}Dir(\boldsymbol{\alpha})$ with the understanding that $\boldsymbol{Y\in
	\Omega}_{p}$ and $\boldsymbol{Y}^{^{\prime}}\in\mathcal{S}_{p+1}$ where%

\[
\boldsymbol{\Omega}_{p}=\left\{  (y_{1},\ldots,y_{p})\in\mathcal{R}^{p}%
:\sum_{i=1}^{p}y_{i}<1,y_{i}>0,i=1,...,p\right\}  ,
\]

\[
\boldsymbol{S}_{p+1}=\left\{  (y_{1},\ldots,y_{p+1})\in\mathcal{R}^{p+1}%
:\sum_{i=1}^{p+1}y_{i}=1,y_{i}>0,i=1,...,p+1\right\}  .
\]

\noindent For any $\boldsymbol{\alpha}$ with \ $\alpha_{i}>0,i=1,...,p+1$ and
$y_{p+1}=1-\sum_{i=1}^{p}y_{i},$ the Dirichlet integral is:%

\begin{equation}%
{\displaystyle\int\limits_{\boldsymbol{\Omega}_{p}}}
{\displaystyle\prod\limits_{i=1}^{p+1}}
y_{i}^{\alpha_{i}-1}d\boldsymbol{y=}%
{\displaystyle\int}
\cdots%
{\displaystyle\int\limits_{\boldsymbol{\Omega}_{p}}}
{\displaystyle\prod\limits_{i=1}^{p+1}}
y_{i}^{\alpha_{i}-1}dy_{1}\cdots dy_{p}=B\left(  \boldsymbol{\alpha}\right)
=\frac{\prod_{i=1}^{p+1}\Gamma(\alpha_{i})}{\Gamma(\alpha_{+})}%
.\label{Multivariate beta function}%
\end{equation}

\noindent \cite{De Groot} and \cite{Kotz} provide detailed discussions on the properties of the Dirichlet distribution.

Assume the baseline distributions to be Gamma$(\theta_{i},\beta_{i}),$
$i=1,...,p,$ with cdfs%

\begin{equation}
G_{i}(x_{i})=\frac{1}{\theta_{i}^{\beta_{i}}\Gamma(\beta_{i})}\int_{0}%
^{x_{i}}e^{-\frac{t}{\theta_{i}}}t^{\beta_{i}-1}dt,\quad\theta_{i},\beta
_{i}>0,\;i=1,\cdots,p,\label{Gamma cdf}%
\end{equation}

\noindent for this chapter. The gamma distribution, which belongs to the
exponential class, is a flexible distribution model with shape parameter
$\beta$, that may offer a good fit to some sets of data.

\section{Recipe}
The construction methodology for the proposed model is as follows:

\begin{itemize}
	\pagebreak
	\item Builder 4:
\end{itemize}%

\begin{equation}
H(x_{1},...,x_{p})=%
{\displaystyle\int_{0}^{G_{1}(x_{1})}}
\cdots%
{\displaystyle\int_{0}^{G_{p}(x_{p})}}
\dfrac{1}{B(\boldsymbol{\alpha})}y_{1}^{\alpha_{1}-1}\cdots y_{p}^{\alpha
	_{p}-1}\left(  1-\sum_{i=1}^{p}y_{i}\right)  ^{\alpha_{p+1}-1}d\boldsymbol{y}%
\label{Builder 4}%
\end{equation}

\noindent where $G_{i}(\cdot)$, $i=1,\ldots,p$, can be any cdf.

Let the joint pdf of $G_{i}(\cdot)$, $i=1,\ldots,p$, be the
Dirichlet pdf given by%

\begin{equation}%
\begin{array}
[c]{rl}%
&f(G_{1},...,G_{p})\\
= & \dfrac{1}{B(\boldsymbol{\alpha})}G_{1}^{\alpha_{1}%
	-1}(x_{1})\ldots G_{p+1}^{\alpha_{p+1}-1}(x_{p+1}),\text{ \ \ \ \ }%
0<G_{i}\left(  \cdot\right)  <1\text{,}\sum_{i=1}^{p+1}G_{i}=1\\
= & \dfrac{1}{B(\boldsymbol{\alpha})}G_{1}^{\alpha_{1}-1}(x_{1})\ldots
G_{p}^{\alpha_{p}-1}(x_{p})\left(  1-\sum_{i=1}^{p}G_{i}(x_{i})\right)
^{\alpha_{p+1}-1},0<\sum_{i=1}^{p}G_{i}(x_{i})<1,\medskip
\end{array}
\label{Joint generator}%
\end{equation}

\noindent i.e. the Dirichlet combines the marginals $G_{i}(\cdot)$,
$i=1,\ldots,p$, with parameters $\boldsymbol{\alpha}=$ $(\alpha_{1}%
,\cdots,\alpha_{p};\alpha_{p+1})$ for \ $\alpha_{i}>0,i=1,...,p+1.\medskip$

\noindent Then, according to $\left(  \ref{Mother}\right)  $, the joint generated distribution, namely the Dirichlet-Gamma ($DG$) has pdf $\medskip$%

\begin{equation}
h(\boldsymbol{x})=\frac{1}{B(\boldsymbol{\alpha})}\left(  1-\sum_{i=1}%
^{p}G_{i}(x_{i})\right)  ^{\alpha_{p+1}-1}\prod_{i=1}^{p}g_{i}(x_{i}%
)G_{i}^{\alpha_{i}-1}(x_{i}),\label{DG pdf}%
\end{equation}

\noindent for $\mathcal{R}^{p}$, $0<\sum_{i=1}^{p}G_{i}(x_{i})<1$ and the parameters $\alpha_{i},\theta_{i},\beta_{i}$ ,$i=1,\cdots,p,$ are
restricted to take those values for which $\left(  \ref{DG pdf}\right)  $ is non-negative, enote $\left(  \ref{DG pdf}\right)  $ as $\boldsymbol{X\sim}%
DG(\boldsymbol{\alpha,\theta,\beta}). $\medskip

Then, the marginal pdf of $X_{i}$, $i=1,...,p,$ has the form%

\begin{equation}
h_{i}(x_{i})=\frac{1}{B(\alpha_{i},\alpha_{+}-\alpha_{i})}g_{i}(x_{i}%
)G_{i}^{\alpha_{i}-1}(x_{i})\left(  1-G_{i}(x_{i})\right)  ^{\alpha_{+}%
	-\alpha_{i}-1},\label{Marginal pdf of DG}%
\end{equation}

\noindent this is useful for determining the moments of $X_{i}$, $i=1,...,p,$.$\medskip$

Although the baseline cdf 's $G_{i}(\cdot)$ could be presented by several
distributions in this chapter, the case where $g_{i}(\cdot)$ is the pdf $\mathrm{Gamma}(\theta_{i},\beta_{i})$, $i=1,\cdots,p$ is considered.

\section{Properties}
Firstly an expression for the product moments will be derived, followed by the moment generating function (mgf) of the $DG(\boldsymbol{\alpha,\theta,\beta})$ distribution. For
this purpose, the following lemma is derived.%

\begin{lemma} 
	\begin{equation}
	\mathcal{I}(\zeta)=%
	{\displaystyle\int}
	\cdots%
	{\displaystyle\int\limits_{\boldsymbol{\Omega}_{p}}}
	\prod_{i=1}^{p}u_{i}^{\alpha_{i}-1}\left(  1-\sum_{i=1}^{p}u_{i}\right)
	^{\zeta}d\boldsymbol{u}\label{Dirichlet integral y1-yp-1}%
	\end{equation}

	\noindent where $\boldsymbol{u}=(u_{1},\ldots,u_{p})$. Then%
	
	\begin{equation}
	\mathcal{I}(\zeta)=\prod_{i=1}^{p-1}B\left(  \alpha_{i},\sum_{j=i+1}^{p}%
	\alpha_{j}+\zeta+1\right)  B\left(  \alpha_{p},\text{ }\zeta+1\right)
	.\label{Expression}%
	\end{equation}
\end{lemma}

\bigskip

\begin{proof}%
	\[%
	\begin{array}
	[c]{rl}%
	\mathcal{I}(\zeta)= & \int_{\boldsymbol{\Omega}_{p}}\prod_{i=1}^{p}%
	u_{i}^{\alpha_{i}-1}\left(  1-\sum_{i=1}^{p}u_{i}\right)  ^{\zeta}\prod
	_{i=1}^{p}du_{i}\bigskip\\
	= & \int_{\boldsymbol{\Omega}_{p}}u_{1}^{\alpha_{1}-1}\prod_{i=2}^{p}%
	u_{i}^{\alpha_{i}-1}(1-u_{1})^{\zeta}\bigskip\\
	& \times\left(  1-\sum_{i=1}^{p}\frac{u_{i}}{1-u_{1}}\right)  ^{\zeta}%
	\prod_{i=1}^{p}du_{i}.\bigskip
	\end{array}
	\]
	
	\noindent Now apply the transformation $v_{i}=\frac{u_{i}}{1-u_{1}}$, for
	$i=2,\cdots,p$, with $J(u_{2},\cdots,u_{p}\rightarrow v_{2},\cdots
	,v_{p})=(1-u_{1})^{p-1}$ to obtain%
	
	\[
	u_{2}=v_{2}(1-u_{1}),\quad\prod_{i=2}^{p}u_{i}^{\alpha_{i}-1}=(1-u_{1}%
	)^{\sum_{i=2}^{p}\alpha_{i}-(p-1)}\prod_{i=2}^{p}v_{i}^{\alpha_{i}-1}.
	\]

	\noindent Hence this results in%
	
	\begin{align*}
	\mathcal{I}(\zeta)  & =\int_{0}^{1}u_{1}^{\alpha_{1}-1}(1-u_{1})^{\zeta
		+\sum_{i=2}^{p}\alpha_{i}}du_{1}\bigskip\\
	& \times\int_{\boldsymbol{\Omega}^{p-1}}\prod_{i=2}^{p}v_{i}^{\alpha_{i}%
		-1}\left(  1-\sum_{i=2}^{p}v_{i}\right)  ^{\zeta}\prod_{i=2}^{p}dv_{i}%
	\bigskip\\
	& =B\left(  \alpha_{1},\sum_{i=2}^{p}\alpha_{i}+\zeta+1\right)  \bigskip\\
	& \times\int_{\boldsymbol{\Omega}^{p-1}}v_{2}^{\alpha_{2}-1}\prod_{i=3}%
	^{p}v_{i}^{\alpha_{i}-1}(1-v_{2})^{\zeta}\left(  1-\sum_{i=3}^{p}\frac{v_{i}%
	}{1-v_{2}}\right)  ^{\zeta}\prod_{i=2}^{p}dv_{i}.
	\end{align*}
	$\bigskip$
	
	\noindent At this stage making the transformation $w_{i}=\frac{v_{i}}{1-v_{2}%
	}$ once more, for $i=3,\cdots,p$, with Jacobian equal to $(1-v_{2})^{p-2}$, it follows that $\bigskip$
	\begin{align*}
	\mathcal{I}(\zeta)  & =B\left(  \alpha_{1},\sum_{i=2}^{p}\alpha_{i}%
	+\zeta+1\right)  B\left(  \alpha_{2},\sum_{i=3}^{p}\alpha_{i}+\zeta+1\right)
	\bigskip\\
	& \times\int_{\boldsymbol{\Omega}^{p-2}}\prod_{i=3}^{p}w_{i}^{\alpha_{i}%
		-1}\left(  1-\sum_{i=3}^{p}w_{i}\right)  ^{\zeta}\prod_{i=3}^{p}%
	dw_{i}.\bigskip
	\end{align*}
	Continuing this procedure, finally yields$\left(  \ref{Expression}\right)  $.
\end{proof}

The following result for the product moment is stated, assuming the pdf $\left(  \ref{DG pdf}\right)  ,$ holds.

\begin{theorem}
	Let $n_{i}$, $i=1,\ldots,n_{p}$ are positive integer values. Then, the product	moments of $\boldsymbol{X\sim}DG(\boldsymbol{\alpha,\theta,\beta})$ admit the
	following explicit form
	\begin{align*}
	\mathcal{E}  & =E\left[  \prod_{i=1}^{p}X_{i}^{n_{i}}\right]  =\left(
	\prod_{i=1}^{p}\frac{\theta_{i}^{n_{i}}\Gamma(n_{i}+\beta_{i})}{\Gamma
		(\beta_{i})}\right)  \\
	& \times\left(  \prod_{i=1}^{p-1}B\left(  \alpha_{i},\sum_{j=i+1}^{p}%
	\alpha_{j}+\frac{\alpha_{p+1}-1}{p}+1\right)  \right)  B\left(  \alpha
	_{p},\frac{\alpha_{p+1}-1}{p}+1\right)  .
	\end{align*}
	$\bigskip$
\end{theorem}

\noindent\textbf{Proof:} From $\left(  \ref{DG pdf}\right)  $, for
$\boldsymbol{X}=(X_{1},\ldots,X_{p})$, it follows that$\bigskip$%

\[%
\begin{array}
[c]{rl}%
\mathcal{E}= &
{\displaystyle\int\limits_{\mathcal{R}^{p},\sum_{j=1}^{p}G_{j}(x_{j})<1}}
\dfrac{1}{B(\boldsymbol{\alpha})}\prod_{i=1}^{p}x_{i}^{n_{i}}\left(
1-\sum_{j=1}^{p}G_{j}(x_{j})\right)  ^{\alpha_{p+1}-1}\prod_{j=1}^{p}%
g_{j}(x_{j})G_{j}^{\alpha_{j}-1}(x_{j})d\boldsymbol{x}\bigskip\\
= & \dfrac{1}{B(\boldsymbol{\alpha})}%
{\displaystyle\int\limits_{\mathcal{R}^{p},\sum_{j=1}^{p}G_{j}(x_{j})<1}}
\prod_{i=1}^{p}\dfrac{\theta_{i}^{n_{i}}\Gamma(n_{i}+\beta_{i})}{\Gamma
	(\beta_{i})}G_{i}^{\alpha_{i}-1}(x_{i})\left(  1-\sum_{i=1}^{p}G_{i}%
(x_{i})\right)  ^{\frac{\left(  \alpha_{p+1}-1\right)  }{p}}\bigskip\\
& \times\dfrac{1}{\theta_{i}^{n_{i}+\beta_{i}}\Gamma(n_{i}+\beta_{i}%
	)}e^{-\frac{x_{i}}{\theta_{i}}}x_{i}^{n_{i}+\beta_{i}-1}d\boldsymbol{x}%
\bigskip\\
= & E\left\{  \prod_{i=1}^{p}\dfrac{\theta_{i}^{n_{i}}\Gamma(n_{i}+\beta_{i}%
	)}{\Gamma(\beta_{i})}G_{i}^{\alpha_{i}-1}(V_{i})\left(  1-\sum_{i=1}^{p}%
G_{i}(V_{i})\right)  ^{\frac{\left(  \alpha_{p+1}-1\right)  }{p}}\right\}
\bigskip
\end{array}
\]
$\bigskip$

\noindent where $V_{i}\sim Gamma(\theta_{i},n_{i}+\beta_{i})$. Using the fact
that $G_{i}(V_{i})\equiv U_{i}\sim U(0,1)$, it follows that%

\[%
\begin{array}
[c]{cc}%
\mathcal{E}= & \prod_{i=1}^{p}\dfrac{\theta_{i}^{n_{i}}\Gamma(n_{i}+\beta
	_{i})}{\Gamma(\beta_{i})}E\left\{  \prod_{i=1}^{p}U_{i}^{\alpha_{i}-1}\left(
1-\sum_{i=1}^{p}U_{i}\right)  ^{\frac{\left(  \alpha_{p+1}-1\right)  }{p}%
}\right\}  \bigskip\\
= & \prod_{i=1}^{p}\dfrac{\theta_{i}^{n_{i}}\Gamma(n_{i}+\beta_{i})}%
{\Gamma(\beta_{i})}%
{\displaystyle\int\limits_{\boldsymbol{\Omega}_{p}}}
\prod_{i=1}^{p}u_{i}^{\alpha_{i}-1}\left(  1-\sum_{i=1}^{p}u_{i}\right)
^{\frac{\left(  \alpha_{p+1}-1\right)  }{p}}d\boldsymbol{u.}\bigskip
\end{array}
\]

\noindent The theorem is completed by applying the Lemma for $\mathcal{I}\left(
\frac{\alpha_{p+1}-1}{p}\right)  $.\hfill$\ \rule{0.5em}{0.5em}$

\begin{theorem}
	The moment generating function (mgf) of $\boldsymbol{X\sim
	}DG(\boldsymbol{\alpha,\theta,\beta})$ is given by
	\begin{align*}
	M_{\boldsymbol{X}}(\boldsymbol{t})  & =\frac{1}{B(\boldsymbol{\alpha})}%
	\sum_{m=0}^{\infty}\frac{1}{m!}\sum_{n_{1}+n_{2}+\cdots+n_{p}=m}\frac
	{m!}{n_{1}!n_{2}!\cdots n_{p}!}\prod_{i=1}^{p}(t_{i})^{n_{i}}\bigskip\\
	& \times\left(  \prod_{i=1}^{p}\frac{\theta_{i}^{n_{i}}\Gamma(n_{i}+\beta
		_{i})}{\Gamma(\beta_{i})}\right)  \prod_{i=1}^{p-1}B\left(  \alpha_{i}%
	,\sum_{j=i+1}^{p}\alpha_{j}+\frac{\alpha_{p+1}-1}{p}+1\right)  B\left(
	\alpha_{p},\frac{\alpha_{p+1}-1}{p}+1\right)  \bigskip
	\end{align*}
	where, $\boldsymbol{t}=(t_{1},\ldots,t_{p})$, $\boldsymbol{x}=(x_{1}%
	,\ldots,x_{p})$ $\boldsymbol{\theta}=(\theta_{1},\cdots,\theta_{p})$ and
	$\boldsymbol{\beta}=(\beta_{1},\cdots,\beta_{p})$.
\end{theorem}

\noindent\textbf{Proof:} 

It follows that%

\[%
\begin{array}
[c]{rl}%
M_{\boldsymbol{X}}(\boldsymbol{t})= & E\left[  e^{\boldsymbol{tX}^{\top}%
}\right]  \bigskip\\
= &
{\displaystyle\int\limits_{\mathcal{R}^{p},\sum_{j=1}^{p}G_{j}(x_{j})<1}}
e^{\boldsymbol{tx}^{\top}}h(\boldsymbol{x})d\boldsymbol{x}\bigskip\\
= &
{\displaystyle\int\limits_{\mathcal{R}^{p},\sum_{j=1}^{p}G_{j}(x_{j})<1}}
\sum_{m=0}^{\infty}\frac{1}{m!}(\boldsymbol{tx}^{\top})^{m}h(\boldsymbol{x}%
)d\boldsymbol{x}\bigskip\\
= &
{\displaystyle\int\limits_{\mathcal{R}^{p},\sum_{j=1}^{p}G_{j}(x_{j})<1}}
\sum_{m=0}^{\infty}\frac{1}{m!}\sum_{n_{1}+n_{2}+\cdots+n_{p}=m}\frac
{m!}{n_{1}!n_{2}!\cdots n_{p}!}\prod_{i=1}^{p}(t_{i}x_{i})^{n_{i}%
}h(\boldsymbol{x})d\boldsymbol{x}\bigskip\\
= & \frac{1}{B(\boldsymbol{\alpha})}\sum_{m=0}^{\infty}\frac{1}{m!}\sum
_{n_{1}+n_{2}+\cdots+n_{p}=m}\frac{m!}{n_{1}!n_{2}!\cdots n_{p}!}\prod
_{i=1}^{p}(t_{i})^{n_{i}}E\left[  \prod_{i=1}^{p}X_{i}^{n_{i}}\right]
\bigskip\\
& \text{where }\top\ \text{denotes transpose of vector.}%
\end{array}
\]

The result follows by Theorem 1.\hfill$\ \rule{0.5em}{0.5em}\bigskip$

$\bigskip$

\section{The proof of the pudding is...}

The basic construction of the $DG\left(  \boldsymbol{\alpha,\theta,\beta
}\right)  $ model \ entails embedding the cdf of a gamma distribution
within the pdf of the Dirichlet distribution, that acts as a generator. The exact generation procedure for the Dirichlet-Gamma random variates is given as Algorithm 1 follows:$\bigskip$%

\begin{tabular}
	[c]{cc}\hline
	& \textbf{Algorithm 1}\\\hline
	\multicolumn{1}{l}{{\small Step 1:}} & \multicolumn{1}{l}{{\small Generate
			independent gamma random variables }$W_{1},W_{2},\ldots,W_{p+1}${\small \ }}\\
	\multicolumn{1}{l}{{\small \ \ \ \ \ \ \ \ \ \ \ \ }} &
	\multicolumn{1}{l}{{\small where }$W_{i}${\small \ }$\sim Gamma\left(
		\alpha_{i},1\right)  ${\small \ for }$\alpha_{i}>0,${\small \ }$i=1,2,\ldots
		,p+1;$}\\
	\multicolumn{1}{l}{{\small Step 2:}} & \multicolumn{1}{l}{{\small Set }%
		$Y_{i}=\frac{W_{i}}{%
			{\textstyle\sum\limits_{j=1}^{p}}
			W_{j}}${\small \ for }$i=1,2,\ldots,p;$}\\
	\multicolumn{1}{l}{{\small Step 3: \ }} & \multicolumn{1}{l}{{\small Return
		}$(Y_{1},Y_{2},\ldots,Y_{p})${\small \ and let }$(Y_{1},Y_{2},\cdots
		,Y_{p})\equiv(G_{1}\left(  x_{1}\right)  ,G_{2}\left(  x_{2}\right)
		,\ldots,G_{p}\left(  x_{p}\right)  )$}\\
	\multicolumn{1}{l}{{\small \ \ \ \ \ \ \ \ \ \ \ \ }} &
	\multicolumn{1}{l}{{\small with }$\sum_{i=1}^{p}G_{i}(x_{i})<1,${\small where
		}$G_{i}\left(  x_{i}\right)  ${\small \ is the cdf of the gamma distribution;}%
	}\\
	\multicolumn{1}{l}{{\small Step 4:}} & \multicolumn{1}{l}{{\small Set }%
		$X_{i}=G_{i}^{-1}\left(  y_{i}\right)  ${\small \ for }$i=1,2,\ldots,p;$}\\
	\multicolumn{1}{l}{{\small Step 5:}} & \multicolumn{1}{l}{{\small Return
		}$(X_{1},X_{2},\cdots,X_{p})${\small \ where }$X\sim DG\left(
		\boldsymbol{\alpha,\theta,\beta}\right)  ${\small \ for}}\\
	\multicolumn{1}{l}{{\small \ \ \ \ \ \ \ \ \ \ \ }$.$} &
	\multicolumn{1}{l}{{\small parameters }$\alpha_{i},\theta_{j},\beta_{j}%
		>0,${\small \ }$i=1,2,\ldots,p+1;j=1,2,\ldots,p.$}%
\end{tabular}

\subsection{Model presentation}\label{sec3}
In Figures 1-6, various pdfs and contour plots of $\left(  \ref{DG pdf}\right)
$ for different values of $\left(  \boldsymbol{\alpha,\theta,\beta}\right)  $
are provided. A 1000 simulated Dirichlet-Gamma values accompany the graphs.

\subsection{Simulation study 1}
Suppose $N$ vector observations $\boldsymbol{X}_{1},\ldots,\boldsymbol{X}_{N}$
of dimension $\left(  p-1\right)  \times1$ are drawn independently and
identically from the $DG(\boldsymbol{\alpha,\theta,\beta})$ distribution.
Therefore, the log-likelihood of $\boldsymbol{\psi}=(\boldsymbol{\alpha
	,\theta,\beta})$ based on the observed data $\left\{  \boldsymbol{X}%
_{i}\right\}  _{i=1}^{N}$ \ from (\ref{DG pdf}) is%

\[
\mathit{l}\left(  \boldsymbol{\psi}\right)  =%
{\displaystyle\sum\limits_{i=1}^{N}}
\log h(\boldsymbol{x;\psi}).
\]

The above simulation Algorithm 1 is used to generate samples of size $100$, $500$ and $1000$. Using $1000$ trials for each group of fixed parameters, 1000 ML estimates of the model parameters (using the optim procedure in R software) is obtained. 

To investigate the estimation accuracies, calculate the mean, bias and mean square error (MSE), defined as%

\[
\text{Bias}=\frac{1}{1000}%
{\displaystyle\sum\limits_{k=1}^{1000}}
\widehat{\psi}_{k}-\psi_{true}\ \ \ \text{and\ \ \ \ MSE}=\frac{1}{1000}%
{\displaystyle\sum\limits_{k=1}^{1000}}
\left(  \widehat{\psi}_{k}-\psi_{true}\right)  ^{2},
\]

\noindent are calculated, where $\widehat{\psi}_{k}$ denotes the ML estimate of $\psi_{true}$
(a specific parameter) at the $k^{th}$ replication. The detailed numerical
results are reported in Table 1-3.

For a large sample size the asymptotic distribution of the ML estimates can be used to construct asymptotic confidence intervals. The asymptotic distribution of the ML estimate of $\psi$ is%

\[
\dfrac{\widehat{\psi}-\psi}{\sqrt{Var\left(  \widehat{\psi}\right)  }}\sim
N\left(  0,1\right)  .
\]

\noindent Confidence intervals (CI) for the model parameters by implementing the parametric bootstrap method are also provided.  Tables 1-3 reflect also the coverage probabilities (CP) and average lengths of the intervals based on these two methods.

\bigskip
\noindent
\textbf{Table 1} Results for $n=100$ and $\boldsymbol{\psi=}\left(  \alpha
_{1},\alpha_{2},\alpha_{3},\beta_{1},\beta_{2},\theta_{1},\theta_{2}\right)
=(2,2,3,1.5,2.8,1.1,1.2).$%

\noindent
\begin{tabular}
	[c]{llllllll}\hline
	$n=100$ & $\widehat{\alpha}_{1}$ & $\widehat{\alpha}_{2}$ & $\widehat{\alpha
	}_{3}$ & $\widehat{\beta}_{1}$ & $\widehat{\beta}_{2}$ & $\widehat{\theta}%
	_{1}$ & $\widehat{\theta}_{2}$\\\hline
	\ Mean & $2.158$ & $2.123$ & $2.861$ & $1.604$ & $2.975$ & $1.222$ & $1.307$\\
	\ Bias & $0.158$ & $0.123$ & $-0.139$ & $0.104$ & $0.175$ & $0.122$ &
	$0.107$\\
	\ MSE & $0.861$ & $0.711$ & $0.688$ & $0.297$ & $0.723$ & $0.144$ & $0.137$\\
	$\text{CP asymptotic CI}$ & $0.945$ & $0.943$ & $0.961$ & $0.949$ & $0.946$ &
	$0.927$ & $0.945$\\
	CP bootstrapped CI & $0.967$ & $0.964$ & $0.966$ & $0.965$ & $0.971$ & $0.965
	$ & $0.968$\\
	\ Length of asymptotic CI & $3.585$ & $3.273$ & $3.206$ & $2.098$ & $3.263$ &
	$1.406$ & $1.391$\\
	\ Length $\text{of bootstrapped CI}$ & $2.876$ & $2.604$ & $5.897$ & $2.839$ &
	$4.894$ & $1.975$ & $1.762$\\\hline
\end{tabular}
\ $\medskip$

\bigskip
\noindent
\textbf{Table 2}: Results for $n=500$ and $\boldsymbol{\psi=}\left(
\alpha_{1},\alpha_{2},\alpha_{3},\beta_{1},\beta_{2},\theta_{1},\theta
_{2}\right)  =(2,2,3,1.5,2.8,1.1,1.2).$%

\noindent
\begin{tabular}
	[c]{llllllll}\hline
	$n=500$ & $\widehat{\alpha}_{1}$ & $\widehat{\alpha}_{2}$ & $\widehat{\alpha
	}_{3}$ & $\widehat{\beta}_{1}$ & $\widehat{\beta}_{2}$ & $\widehat{\theta}%
	_{1}$ & $\widehat{\theta}_{2}$\\\hline
	\ Mean & $2.018$ & $2.025$ & $2.928$ & $1.535$ & $2.851$ & $1.138$ & $1.234$\\
	\ Bias & $0.018$ & $0.025$ & $-0.072$ & $0.035$ & $0.051$ & $0.038$ &
	$0.034$\\
	\ MSE & $0.182$ & $0.161$ & $0.111$ & $0.054$ & $0.156$ & $0.028$ & $0.027$\\
	$\text{CP asymptotic CI}$ & $0.946$ & $0.939$ & $0.948$ & $0.950$ & $0.938$ &
	$0.940$ & $0.937$\\
	$\text{CP bootstrapped CI}$ & $0.974$ & $0.974$ & $0.975$ & $0.975$ & $0.975$
	& $0.975$ & $0.974$\\
	$\ $Length$\text{ of asymptotic CI}$ & $1.670$ & $1.571$ & $1.278$ & $0.903$ &
	$1.537$ & $0.635$ & $0.627$\\
	$\ $Length $\text{of bootstrapped CI}$ & $2.124$ & $1.951$ & $2.848$ & $1.478
	$ & $2.755$ & $1.207$ & $1.341$\\\hline
\end{tabular}
\ $\medskip$

\pagebreak
\noindent
\textbf{Table 3}: Results for $n=1000$ and $\boldsymbol{\psi=}\left(
\alpha_{1},\alpha_{2},\alpha_{3},\beta_{1},\beta_{2},\theta_{1},\theta
_{2}\right)  =(2,2,3,1.5,2.8,1.1,1.2).$%

\noindent
\begin{tabular}
	[c]{llllllll}\hline
	$n=1000$ & $\widehat{\alpha}_{1}$ & $\widehat{\alpha}_{2}$ & $\widehat{\alpha
	}_{3}$ & $\widehat{\beta}_{1}$ & $\widehat{\beta}_{2}$ & $\widehat{\theta}%
	_{1}$ & $\widehat{\theta}_{2}$\\\hline
	Mean & $1.999$ & $2.005$ & $2.963$ & $1.526$ & $2.837$ & $1.125$ & $1.224$\\
	Bias & $-0.001$ & $0.005$ & $-0.038$ & $0.026$ & $0.037$ & $0.025$ & $0.024$\\
	MSE & $0.099$ & $0.082$ & $0.059$ & $0.032$ & $0.083$ & $0.016$ & $0.015$\\
	$\text{CP asymptotic CI}$ & $0.946$ & $0.944$ & $0.942$ & $0.937$ & $0.938$ &
	$0.927$ & $0.935$\\
	$\text{CP bootstrapped CI}$ & $0.975$ & $0.975$ & $0.975$ & $0.975$ & $0.975$
	& $0.975$ & $0.975$\\
	Length$\text{ of asymptotic CI}$ & $1.237$ & $1.125$ & $0.941$ & $0.690$ &
	$1.118$ & $0.488$ & $0.478$\\
	Length$\text{ of bootstrapped CI}$ & $2.010$ & $2.020$ & $2.847$ & $1.456$ &
	$2.754$ & $1.237$ & $1.270$\\\hline
\end{tabular}
\bigskip

It can be observed \ that the bias and MSE of the $\boldsymbol{DG}\left(
\boldsymbol{\alpha,\theta,\beta}\right)  $ distribution tend to decrease
toward zero by increasing sample size ($n$), showing empirically the
consistency of the ML estimates. The MSE of the estimates of $\widehat
{\boldsymbol{\beta}}$ is higher than  {$\widehat{\boldsymbol{\theta}}$}, as one would expect from the shape parameter of the gamma baselines. As the sample size changes from 100 to 1000, the average length of confidence intervals do decrease.

\subsection{Simulation study 2}

A model testing technique, referred to in this chapter as the empirical
estimator of the cdf of a multivariate distribution, is proposed in analysing the performances of the two competing models, namely the Dirichlet ($D$) and Dirichlet-Gamma ($DG$). The technique compares the empirical cdfs of the observed and simulated datasets. The following steps ( Algorithm 2) are taken in order to assess the competence of the models.%

The advantage of this technique, is that one can also use the empirical cdfs to rank the simulated data.  Ranking data makes it possible to calculate more accurate distances between the observed data points and the simulated points.  Figure 7 illustrates an observed dataset (in black) and simulated points from the simulated artificial datasets Dirichlet (in blue) and the Dirichlet-Gamma (in red). \ The challenge lies in choosing the correct simulated point to calculate the distances.  The solution that is proposed in this chapter is to rank the simulated data from the two competing models according to their calculated empirical cdfs respectively. The distances (as shown with the arrows) between the observed (in black) and the simulated data points can be more accurately calculated based on the quantile positions.

\begin{tabular}
	[c]{cc}\hline
	& \textbf{Algorithm 2}\\\hline
	\multicolumn{1}{l}{{\small Step 1:}} & \multicolumn{1}{l}{{\small From the
			observed dataset }$x_{n\times p}${\small , calculate the empirical cdf }}\\
	\multicolumn{1}{l}{{\small \ }} & \multicolumn{1}{l}{$\widehat{F\left(
			\mathbf{x}\right)  }=P\left(  X_{1}\leq x_{1},X_{2}\leq x_{2},\cdots,X_{p}\leq
		x_{p}\right)  =\frac{1}{n}%
		{\textstyle\sum\limits_{i=1}^{p}}
		I\left(  x_{i}\leq x\right)  ,${\small \ }}\\
	& \multicolumn{1}{l}{{\small where }$I\left(  \cdot\right)  ${\small \ is the
			indicator function;}}\\
	\multicolumn{1}{l}{{\small Step 2:}} & \multicolumn{1}{l}{{\small Obtain the
			parameter estimates for the two competing models, }${\small D}${\small \ and
			\ }${\small DG}$}\\
	& \multicolumn{1}{l}{{\small distributions and simulate artificial datasets;}%
	}\\
	& \multicolumn{1}{l}{$x_{D}^{\ast}=\left(  x_{1}^{\ast},x_{2}^{\ast
		},\cdots,x_{p}^{\ast}\right)  ${\small \ and }$x_{DG}^{\ast
		}=\left(  x_{1}^{\ast},x_{2}^{\ast},\cdots,x_{p}^{\ast}\right)  ${\small \ of sizes }$d>n.$}\\
	\multicolumn{1}{l}{{\small Step 3: \ }} & \multicolumn{1}{l}{{\small Calculate
			the empirical cdfs for each simulated artificial dataset}}\\
	\multicolumn{1}{l}{} & \multicolumn{1}{l}{$\widehat{F\left(  \mathbf{x}^{\ast
			}\right)  }=P\left(  X_{1}^{\ast}\leq x_{1},X_{2}^{\ast}\leq x_{2}%
		,\cdots,X_{p}^{\ast}\leq x_{p}\right)  =\frac{1}{d}%
		{\textstyle\sum\limits_{i=1}^{p}}
		I\left(  x_{i}^{\ast}\leq x\right)  ${\small ;}}\\
	\multicolumn{1}{l}{{\small Step 4:}} & \multicolumn{1}{l}{{\small Repeat step
			2 - 3 }$m${\small \ times, and for each simulation, compute
			Kolmogorov-Smirnov (KS) }}\\
	\multicolumn{1}{l}{} & \multicolumn{1}{l}{{\small distances between the
			empirical cdf (as computed in step 1) and the empirical cdfs}}\\
	\multicolumn{1}{l}{} & \multicolumn{1}{l}{{\small of the competing models (as
			computed in step 3) where KS measure is defined in this case as}}\\
	& \multicolumn{1}{l}{$KS=\max\left\vert \widehat{F\left(  \mathbf{x}^{\ast
			}\right)  }-\widehat{F\left(  \mathbf{x}\right)  }\right\vert $}\\
	\multicolumn{1}{l}{{\small Step 5:}} & \multicolumn{1}{l}{{\small Compute the
			average KS distances over the }$m${\small \ simulated artificial datasets;}%
	}\\
	\multicolumn{1}{l}{{\small Step 6:}} & \multicolumn{1}{l}{{\small Compare the
			KS distances of the }${\small DG}${\small \ to the KS distance}}\\
	\multicolumn{1}{l}{} & \multicolumn{1}{l}{{\small of the }${\small D}%
		${\small \ in terms of the ratio }$\frac{KS\text{ of }DG\text{ }}{KS\text{ of
			}D\text{ }}.$}\\
	&
\end{tabular}

In this chapter  for the implementation of this technique, the focus is on the ratio of the KS distances between the two competing models. 
\noindent To test this model testing technique, generate a "observed" dataset from a Dirichlet distribution and analyse the performance of the
Dirichlet-Gamma through the steps. \ Since the KS distances vary
from simulation to simulation, samples of sizes $d=100,1000,10000$ are
generated from the obtained parameter estimates for Dirichlet and
Dirichlet-Gamma from the observed , where KS distances are calculated for each simulated dataset group. \ $\bigskip$

\noindent i. Generate an artificial dataset from the Dirichlet distribution with parameters $\left(  \alpha_{1},\alpha_{2},\alpha_{3}\right)  =\left(
2,2,3\right)  $ and assume it as the observed data;

\bigskip
\noindent ii. Using this observed dataset, obtain parameter estimates for the Dirichlet and Dirichlet-Gamma distributions;
\bigskip

\noindent iii. From the obtained parameter estimates simulate datasets
of sizes $d=100,1000,\\10000$. Calculate the empirical cdfs for each simulation,  as seen in step 3 of Algorithm 2;
\bigskip


\noindent iv. Calculate the KS distances between the empirical cdf and the cdfs of the two competing models, for each group;
\bigskip

\noindent v. Repeat steps (iii.-iv.) a $100$ times and compute the average KS distance for the two models.
\bigskip

\noindent vi. Represent the KS distance of the Dirichlet-Gamma and Dirichlet as
a ratio $\frac{KS\text{ } of DG\text{ }}{KS\text{ } of D\text{ }}$ for each simulated group of $d=100,1000,10000$.

\medskip
\noindent 
It is observed in Figure 8 that the Dirichlet-Gamma distribution is flexible enough to model Dirichlet distributed variables. The KS distance of the Dirichlet-Gamma is seen to be smaller for all simulated groups.
\subsection{Simulation study 3}
\noindent A further simulation study is carried out to illustrate the
flexibility of the Dirichlet-Gamma when outliers are present within a dataset. Suppose that two non-Dirichlet artificial compositional datasets, where outliers are present, are generated, using Algorithm 3.\bigskip%

\begin{tabular}
	[c]{cc}\hline
	& \textbf{Algorithm 3}\\\hline
	\multicolumn{1}{l}{{\small Step 1:}} & \multicolumn{1}{l}{{\small Generate
		}$n${\small \ random variates }$W_{i}${\small \ }$\symbol{126}Weibull\left(
		k_{i},\lambda_{i}\right)  ${\small \ for }$i=1,2,3.$}\\
	\multicolumn{1}{l}{{\small Step 2:}} & \multicolumn{1}{l}{{\small Define
			random variables }$Y=\left(  Y_{1},Y_{2},Y_{3}\right)  ,${\small \ where
		}$Y_{i}=\frac{W_{i}}{%
			{\textstyle\sum\limits_{i=1}^{3}}
			W_{i}}${\small , }$i=1,2,3${\small .}}\\
	& \multicolumn{1}{l}{{\small and generate artificial dataset }$y=(y_{1}%
		,y_{2},y_{3})$}%
\end{tabular}
\bigskip

\noindent The construction of random variables $Y_{1},Y_{2},Y_{3}$ yields a compositional dataset with a negative correlation. \ The initial values for the Dirichlet and Dirichlet-Gamma used in the R package \textit{optim} are obtained through a grid search. \ Figures 9 and 10 illustrates the flexibility of the Dirichlet-Gamma over outliers.\medskip



\subsection{Real data analysis}

To investigate the performance of the Dirichlet-Gamma distribution with
respect to the Dirichlet distribution, different goodness-of-fit measures will be used to evaluate the models as candidates for the different datasets, namely the Q-Q plot, the Akaike information criterion (AIC, \cite{Akaike}) and the Bayesian information criterion (BIC, \cite{Schwarz}), with the last 2 measures defined as%

\[
\text{AIC}=2m-2\mathit{l}_{max}\ \ \text{and\ \ BIC}=m\log N-2\mathit{l}%
_{max},
\]

\noindent where $m$ is the number of free parameters and $\mathit{l}_{max}$ is the maximized log-likelihood value. Models with lower values of AIC\ and BIC are considered more preferable.\bigskip

\subsubsection{\textbf{EXAMPLE\ 1-Pekin ducklings dataset}}

As first illustration, the Serum-protein data of white Pekin ducklings are
considered (see \cite{Mosimann (1962)}). To illustrate the performance of the Dirichlet-Gamma model with respect to extreme outlying observations, observation 20 of the dataset was perturbated. The blood serum proportions
(pre-albumin, albumin and globulin) in 3-week-old Pekin ducklings were
reported with correlation matrix:%

\[
\left[
\begin{array}
[c]{ccc}%
1 & -0.108 & -0.557\\
-0.108 & 1 & -0.766\\
-0.557 & -0.766 & 1
\end{array}
\right]  .
\]

Using randomly chosen initial parameter values $(\alpha_{1},\alpha_{2}%
,\alpha_{3})=(6.856,2.392,1)$ and $(\alpha_{1},\alpha_{2},\alpha_{3},\beta_{1}%
,\theta_{1},\beta_{2},\theta_{2})=(2.016,2.757,3.318,0.559,0.826,1.569,1.876)$ to obtain the ML estimates of the Dirichlet and Dirichlet-Gamma respectively with the optim package in R.\ The simulated Dirichlet and Dirichlet-gamma random variates are obtained using the ML estimates. Figure 11 shows the Q-Q plots on distances to origin of observed and Dirichlet simulated data.

Figure 12 shows the observed data (black dots) versus simulated data from the Dirichlet distribution (blue dots), accompanied by a contour plot. It is clear that the Dirichlet distribution does not cover all the data points well. Similarly, the red dots show the simulated Dirichlet - Gamma values with a contour plot (second row on Figure 12). The results presented in Figure 12, illustrates that the Dirichlet-Gamma distribution provides a dataset closer to the observed data compared to the Dirichlet distribution. The Dirichlet-Gamma covers the outlier while the Dirichlet model could not detect it. Table 3 shows a summary of the ML fittings (note Log-likelihood is indicated as
$\mathit{ll}$ in the tables).



\noindent
\textbf{Table 4:} Parameter estimates and the performance summary for the
Pekin duckling dataset%

\noindent
\begin{tabular}
	[c]{lllllllllll}\hline
	Model & \multicolumn{7}{l}{ML estimates} & \multicolumn{3}{l}{}\\\hline
	& $\widehat{\alpha}_{1}$ & $\widehat{\alpha}_{2}$ & $\widehat{\alpha}_{3}$ &
	$\widehat{\beta}_{1}$ & $\widehat{\beta}_{2}$ & $\widehat{\theta}_{1}$ &
	$\widehat{\theta}_{2}$ & $\mathit{ll}$ & AIC & BIC\\\hline
	Dirichlet & 4.786 & 28.798 & 30.653 & n/a & n/a & n/a & n/a & -79.797 &
	165.594 & 169.0015\\\hline
	DG & 2.173 & 2.466 & 13.998 & 0.971 & 1.383 & 6.711 & 8.537 & -63.205 &
	140.409 & 148.358\\\hline
\end{tabular}
\bigskip

Using the model testing technique as described by Algorithm 2, it is observed that the KS distance is smaller in the case of the proposed Dirichlet-Gamma model versus the Dirichlet model (see Figure 13).
\subsubsection{\textbf{EXAMPLE\ 2-White cells dataset}}

Three kind of white cells (granulocytes, lymphocytes, monocytes) found in 30 blood samples are recorded in this dataset. The inputs result in 30 pairs of 3-part compositions of the white cells, where each portion was determined through time-consuming microscopic and automatic image analysis. The correlation matrix is given as%

\[
\left[
\begin{array}
[c]{ccc}%
1 & -0.832 & -0.405\\
-0.832 & 1 & -0.170\\
-0.405 & -0.170 & 1
\end{array}
\right]  .
\]

The Dirichlet and the Dirichlet-Gamma distributions are tested to see if they are suitable contenders of this dataset. Using randomly chosen initial parameter values $(\alpha_{1},\alpha_{2},\alpha_{3})=(1,1,1)$ and $(\alpha
_{1},\alpha_{2},\alpha_{3},\beta_{1},\theta_{1},\beta_{2},\theta
_{2})=(2,3,7,1,1.5,0.5,1)$ in this case. The Q-Q plots, scatter plots and
contour plots are presented in Figures 14 and 15, together with the summary of the results when fitting the Dirichlet and the Dirichlet-Gamma to this dataset. It is observed that the Dirichlet-Gamma outperforms the Dirichlet  model.\bigskip

\noindent
\textbf{Table 5:} Parameter estimates and the performance summary for the
White cells dataset%

\noindent
\begin{tabular}
	[c]{lllllllllll}\hline
	Model & \multicolumn{7}{l}{ML estimates} & \multicolumn{3}{l}{}\\\hline
	& $\widehat{\alpha}_{1}$ & $\widehat{\alpha}_{2}$ & $\widehat{\alpha}_{3}$ &
	$\widehat{\beta}_{1}$ & $\widehat{\beta}_{2}$ & $\widehat{\theta}_{1}$ &
	$\widehat{\theta}_{2}$ & $\mathit{ll}$ & AIC & BIC\\\hline
	Dirichlet & 3.208 & 1.455 & 0.593 & n/a & n/a & n/a & n/a & -51.410 &
	108.820 & 113.023\\\hline
	DG & 25.389 & 4.370 & 1.142 & 0.199 & 0.479 & 0.483 & 0.065 & -30.155 &
	74.310 & 84.118\\\hline
\end{tabular}




\section{Conclusion}

This chapter's broader target was to show that the \textquotedblleft mother technique\textquotedblright\ ( see\ref{beta normal cdf})  can still generate novel progeny. A unique contribution is made by introducing a constructive methodology for families of multivariate distributions through the model $H(\boldsymbol{x})=F(G(\boldsymbol{x}))$ with $\boldsymbol{x}$ a vector; $G(\boldsymbol{x})$ a vector of independent Gamma cdfs referred to as baseline distributions and $F$ a multivariate pdf such as the Dirichlet with negative correlations between variables. Simulation studies and two real life cases are investigated to illustrate the value added of this construction, using several performance measures. A new model testing technique based on the empirical estimator of the cdf, is introduced to evaluate the performance of multivariate models. It flows naturally that instead of the gamma baseline distributions any other family of distributions could be used, similarly a more general structure for the generator could be the Dirichlet-hyper-geometric function type I distribution \cite{Nagar}. To accommodate for positive correlation structure in the data, the authors consider the Dirichlet type III\ distribution (see \cite{Ehlers}) or the Liouville distribution of the second kind (\cite{Gupta and Richards},
\cite{Bouguila}) in a follow-up paper. Note that, in contrast with the
Dirichlet and like the generalized Dirichlet, the covariance can be positive or negative.\ The builder would be of the form:

\begin{itemize}
	\item Builder 5:
\end{itemize}%

\[
H(x_{1},...,x_{p})=%
{\displaystyle\int_{0}^{G_{1}(x_{1})}}
\cdots%
{\displaystyle\int_{0}^{G_{p}(x_{p})}}
C\prod_{i=1}^{p}y_{i}^{\alpha_{i}-1}q\left(  \sum_{i=1}^{p}y_{i}\right)
d\boldsymbol{y}%
\]

\noindent where $G_{i}(\cdot)$, $i=1,\ldots,p$, can be any cdf, $C$\ the
normalizing constant of the pdf of the generator and $q(\cdot)$ a measurable
positive real valued function defined on the interval $\left(  0,1\right)  $
such that $%
{\displaystyle\int\limits_{0}^{1}}
q\left(  \tau\right)  \tau^{s-1}d\tau$ exists for all $s>0.$

\noindent This new approach to construct multivariate distributions expands
the body of knowledge within the distribution theory domain.\bigskip

\section*{Acknowledgements}

We express our sincere thanks to Mehrdad Naderi for many helpful
conversations. This work is based on the research supported in part by the
National Research Foundation of South Africa (Grant ref. CPRR160403161466 nr.
105840 and grant ref. IFR170227223754 nr. 109214). Opinions expressed and
conclusions arrived at are those of the authors and are not necessarily to be attributed to the NRF. The authors would like to thank the reviewers for their valuable contributions.

All figures can be obtained from the authors. 

\section*{Appendix}

Code with comments for this document is available from the corresponding author.


\begin{thebibliography}{99}                                                                                               %
\bibitem {Akaike}Akaike, H (1998). Information theory and an extension of the
maximum likelihood principle. In: Selected papers of Hirotugu Akaike. Springer:199--213.

\bibitem {Alexander}Alexander C, Cordeiro GM, Ortega, EMM (2012). Generalized
beta-generated distributions, Comput Statist Data Anal, 56:1880--1897.

\bibitem {Balakrishnan2003}Balakrishnan, N, Nevzorov, VB (2003). A Primer on
Statistical distributions, John Wiley \& Sons, New York, USA.

\bibitem {Barreto-Souza}Barreto-Souza, W, Santos, AHS and Cordeiro, GM (2010).
The beta generalized-exponential distribution, Journal of Statistical
Computation and Simulation, 80(2): 159-172.

\bibitem {Barndorff-Nielsen}Barndorff-Nielsen,OE, Jorgensen, B (1991). Some
parametric models on the simplex, Journal of Multivariate Analysis 39: 106--116.

\bibitem {Bouguila}Bouguila, N (2011). Count Data Modeling and Classification
Using Finite Mixtures of Distributions, IEEE Transactions On Neural Networks,
22,(2): 186-197.

\bibitem {Connor and Mosimaan}Connor, JR, Mosimann, JE (1969). Concepts of
independence for proportions with a generalization of the Dirichlet
distribution, Journal of the American Statistical Association 64:194--206.

\bibitem {De Groot}De Groot, MH (1970). Optimal Statistical Decisions. McGraw-Hill.

\bibitem {Ehlers}Ehlers, R (2011). Bimatrix variate distributions of Wishart
ratios with application, Unpublished dissertation, University of Pretoria.

\bibitem {Elgarhy 2016}Elgarhy, M, Hassan, A S, Rashed, M (2016).
Garhy-generated family of distributions with application, Mathematical Theory
and Modeling 6 (2):

\bibitem {Epaillard}Epaillard, A, Bouguila, N (2019). Data-free metrics for
Dirichlet and generalized Dirichlet mixture-based HMMs--A practical study,
Pattern Recognition, 8: 207--219.

\bibitem {Eugene}Eugene, N, Lee, C, Famoye, F (2002). Beta-normal distribution
and its applications. Communications in Statistics---Theory and Methods 31(4): 497--512.

\bibitem {Favaro}Favaro, S, Hadjicharalambous, G, Prunster, I (2011). On a
class of distributions on the simplex, Journal of Statistical Planning and
Inference 141: 2987--3004.

\bibitem {Gupta and Richards}Gupta, RD, Richards, D.St.P(1997). Multivariate
Liouville distributions, V, in: N.L. Johnson, N. Balakrishnan (Eds.), Advances
in the Theory and Practice of Statistics: A Volume in Honour of Samuel Kotz,
Wiley, New York, 377--396.

\bibitem {Jones 2004}Jones, MC (2004). Families of distributions arising from
distributions of order statistics, Test 13 (1): 1--43.

\bibitem {Kotz}Kotz, S, Balakrishnan, N, Johnson, NL (2000). Continuous
Multivariate Distributions, Vol. 1, Second Edition, John Wiley \& Sons, New
York, USA

\bibitem {Makgai 2017}Makgai, SL, Bekker, A, Ferreira, JT, Arashi, M (2017).
New results from a beta-Pareto class, South African Statistical Journal, 51: 345-360.

\bibitem {Makgai 2019}Makgai, SL, Visagie, J, Bekker, A, De Waal, D (2019).
Contributions to the class of beta-generated distributions, in preparation to submit.

\bibitem {Mameli}Mameli, V (2015). The Kumaraswamy skew-normal distribution,
Statistics \& Probability Letters,104: 75-81.

\bibitem {Mosimann (1962)}Mosimann, JE (1962). \ On the compound multinomial
distribution, the multivariate $\beta$-distribution, and correlations among
proportions, Biometrika, \textbf{49}: 65-82.

\bibitem {Nagar}Nagar, DK, \textperiodcentered\ Bran-Cardona, Gupta, AK
(2009). Multivariate Generalization of the Hypergeometric Function Type I
Distribution, Acta Appl Math, 105: 111--122

\bibitem {Nassar}Nassar, M, Kumar, Dey, D, Cordeiro, GM, Afify, AZ (2019). The
Marshall Olkin alpha power family of distributions with applications, Journal
of Computational and Applied Mathematics, 351: 41-53.

\bibitem {Nadarajah2006}Nadarajah S, Kotz S (2006). The beta exponential
distribution, Reliab Eng Syst Safe, 91:689--697.

\bibitem {Ng}Ng, KW, Tian, G, Tang, M (2011). Dirichlet and related
distributions, Theory, Methods and Applications, John Wiley \& Sons, New York, USA.

\bibitem {Olkin}Olkin, I, Liu, R(2003). A bivariate beta distribution, Stat.
Probability Lett. 62:407--412.

\bibitem {Ongaro}Ongaro, A, S. Migliorati, S (2013). A generalization of the
Dirichlet distribution, Journal of Multivariate Analysis, 114: 412-426.

\bibitem {Ristic}Risti\'{c}, MM, Popovi\'{c}, BV, Zografos, K, Balakrishnan, N
(2018). Discrimination among bivariate beta-generated distributions,
Statistics, 52:2, 303-320

\bibitem {Samanthi}Samanthi, RGM, Sepanski, J (2017). A bivariate extension of
the beta generated distribution derived from copulas, Communications in
Statistics-Theory and Methods, 0(0): 1-17.

\bibitem {Sarabia 2014}Sarabia, J.M, Prieto, F and V. Jord\'{a}, V (2014).
Bivariate beta-generated distributions with applications to well-being data,
Journal of Statistical Distributions and Applications,1: 1--15.

\bibitem {Schwarz}Schwarz, G (1978). Estimating the dimension of a model,The
Annals of Statistics 6 (2): 461--464.

\bibitem {Thomas}Thomas, S, Jacob, J (2006).. A generalized Dirichlet model,
Statistics and Probability Letters 76: 1761-1767.

\bibitem {Weber}Weber, MD, Leemis, LM, Kinciad, R.K. (2006). Minimum
Kolmogorov-Smirnov test statistic parameter estimates. Journal of Statistical
Computation and Simulation, 76 (3): 195--205.

\bibitem {Zografos}Zografos K, Balakrishnan,N (2009). On families of beta- and
generalized gamma-generated distributions and associated inference. Stat
Methodol, 6: 344--362.
\end{thebibliography}
\end{document}